\documentclass[10pt,a4paper,reqno]{amsart}
\usepackage{amsmath}
\usepackage{amsthm}
\usepackage{amsfonts}
\usepackage{amssymb}

\title{Words and Transcendence}

\author{Michel Waldschmidt}

\address{Universit\'e Pierre et Marie Curie -- Paris 6, UMR 7586, Institut de Math\'ematiques de Jussieu, 175 rue du Chevaleret, F-75013 Paris, France}

\email{miw@math.jussieu.fr}

\dedicatory{Dedicated to the 80th birthday of Professor K.F. Roth}

\newtheorem{theorem}{Theorem}[section]
\newtheorem{conjecture}[theorem]{Conjecture}
\newtheorem{corollary}[theorem]{Corollary}
\newtheorem{proposition}[theorem]{Proposition}
\newtheorem*{hypothesisa}{Hypothesis A}
\theoremstyle{definition}
\newtheorem*{example}{Example}
\theoremstyle{definition}
\newtheorem*{remark}{Remark}
\theoremstyle{definition}
\newtheorem*{acknowledgments}{Acknowledgments}

\def\bfF{\mathbf{F}}
\def\bfN{\mathbf{N}}
\def\bfQ{\mathbf{Q}}
\def\bfZ{\mathbf{Z}}

\def\bfu{\mathbf{u}}
\def\bfv{\mathbf{v}}
\def\bfx{\mathbf{x}}

\def\AAA{\mathcal{A}}
\def\QQQ{\mathcal{Q}}

\renewcommand{\le}{\leqslant}
\renewcommand{\ge}{\geqslant}

\begin{document}

\setcounter{page}{1}

\begin{abstract}
Is it possible to distinguish algebraic from transcendental real numbers by considering the $b$-ary expansion in some base $b\ge2$?
In 1950, \'E. Borel  suggested that the answer is no and that for any real irrational algebraic number $x$ and for any base $g\ge2$, the $g$-ary expansion of $x$ should satisfy some of the laws that are shared by almost all numbers.
For instance, the frequency where a given finite sequence of digits occurs should depend only on the base and on the length of the sequence. 

We are very far from such a goal: there is no explicitly known example of a triple $(g,a,x)$, where
$g\ge3$ is an integer, $a$ a digit in $\{0,\ldots,g-1\}$ and $x$ a real irrational algebraic number, for which one can claim that the digit $a$ occurs infinitely often in the $g$-ary expansion of~$x$.

Hence there is a huge gap between the established theory and the expected state of the art.
However, some progress has been made recently, thanks mainly to clever use of Schmidt's subspace theorem.
We review some of these results.
\end{abstract}

\maketitle

%%%%%%%%%%%%%%%%%%%%
%
% SECTION 1
%
%%%%%%%%%%%%%%%%%%%%

\section{Normal Numbers and Expansion of Fundamental Constants}\label{S:NormalNumbers}

%%%%%%%%%%%%%%%%%%%%
%
% SUBSECTION 1.1
%
%%%%%%%%%%%%%%%%%%%%

\subsection{Borel and Normal Numbers}\label{SS:Borel}

In two papers, the first \cite{waldschmidt-borel09} published in 1909 and the second
\cite{waldschmidt-borel50} in 1950, Borel studied the $g$-ary expansion of real numbers, where
$g\ge2$ is a positive integer.
In his second paper he suggested that this expansion for a real irrational algebraic number should satisfy some of the laws shared by almost all numbers, in the sense of Lebesgue measure.

Let $g\ge2$ be an integer.
Any real number $x$ has a unique expansion
\begin{displaymath}
x=a_{-k}g^{k}+\ldots+a_{-1}g+a_0+a_1g^{-1}+a_2g^{-2}+\ldots,
\end{displaymath}
where $k\ge0$ is an integer and the $a_i$ for $i\ge-k$, namely the digits of $x$ in the expansion in base $g$ of~$x$, belong to the set $\{0,1,\ldots,g-1\}$.
Uniqueness is subject to the condition that the sequence $(a_i)_{i\ge-k}$ is not ultimately constant and equal to $g-1$.
We write this expansion
\begin{displaymath}
x=a_{-k}\ldots a_{-1}a_0.a_1a_2\ldots.
\end{displaymath}

\begin{example}
We have
\begin{displaymath}
\sqrt{2}=1.41421356237309504880168872420\ldots
\end{displaymath}
in base $10$ (\emph{decimal expansion}), whereas
\begin{displaymath}
\sqrt{2}=1.01101010000010011110011001100111111100111011110011\ldots
\end{displaymath}
in base $2$ (\emph{binary expansion}).
\end{example}

The first question in this direction is whether each digit always occurs at least once.

\begin{conjecture}\label{Conj:1}
Let $x$ be an real irrational algebraic number, $g\ge3$ a positive integer and $a$ an integer in the range $0\le a\le g-1$.
Then the digit $a$ occurs at least once in the $g$-ary expansion of~$x$.
\end{conjecture}

For $g=2$ it is plain that each of the two digits $0$ and $1$ occurs infinitely often in the binary expansion of an irrational algebraic number.
The same is true for each of the sequences of two digits $01$ and~$10$.
Conjecture \ref{Conj:1} implies that each of the sequences of digits $00$ and $11$ should also occur infinitely often in such an expansion; apply Conjecture \ref{Conj:1} with $g=4$.
There is no explicitly known example of a triple $(g,a,x)$, where $g\ge3$ is an integer, $a$ a digit in
$\{0,\ldots, g-1\}$ and $x$ a real irrational algebraic number, for which one can claim that the digit $a$ occurs infinitely often in the $g$-ary expansion of~$x$.
Another open problem is to produce an explicit pair $(x,g)$, where $g\ge3$ is an integer and $x$ a real irrational algebraic number, for which we can claim that the number of digits which occur infinitely many times in the $g$-ary expansion of $x$ is at least~$3$.
Even though few results are known and explicit examples are lacking, something is known.
For any $g\ge2$ and any $k\ge1$ there exist real algebraic numbers $x$ such that any sequence of $k$ digits occurs infinitely often in the $g$-ary expansion of~$x$.
However the connection with algebraicity is weak: indeed Mahler proved in 1973 more precisely that for any real irrational number $\alpha$ and any sequence of  $k$ digits in the set $\{0,\ldots,g-1\}$, there exists an integer $m$ for which the $g$-ary expansion of $m\alpha$ contains infinitely many times the given sequence; see \cite[Theorem~M]{waldschmidt-AB05b} and~\cite{waldschmidt-mahler73}.
According to Mahler, the smallest such $m$ is bounded by $g^{2k+1}$.
This estimate has been improved by Berend and Boshernitzan \cite{waldschmidt-BB94} to $2g^{k+1}$ and one cannot get better than $g^k-1$.

If a real number $x$ satisfies Conjecture \ref{Conj:1} for all $g$ and~$a$, then it follows that for
any~$g$, each given sequence of digits occurs infinitely often in the $g$-ary expansion of~$x$.
This is easy to see by considering powers of~$g$.

Borel asked more precise questions on the frequency of occurrences of sequences of binary digits of real irrational algebraic numbers. 
We need to introduce some definitions.

Firstly, a real number $x$ is called \emph{simply normal in base $g$} if each digit occurs with frequency $1/g$ in its $g$-ary expansion.
A very simple example in base $10$ is
\begin{displaymath}
x=0.123456789012345678901234567890\ldots,
\end{displaymath}
where the sequence $1234567890$ is repeated periodically, but this number is rational.
We have
\begin{displaymath}
x
=\frac{1\,234\,567\,890}{9\,999\,999\,999}
=\frac{137\,174\,210}{1\,111\,111\,111}.
\end{displaymath}

Secondly, a real number $x$ is called \emph{normal in base $g$} or \emph{$g$-normal} if it is simply normal in base $g^m$ for all $m\ge1$.
Hence a real number $x$ is normal in base $g$ if and only if, for all $m\ge1$, each sequence of $m$ digits occurs with frequency $1/g^m$ in its $g$-ary expansion.

Finally, a number is called \emph{normal} if it is normal in all bases $g\ge2$.

Borel suggested in 1950 that each real irrational algebraic number should be normal.

\begin{conjecture}[Borel, 1950]\label{Conj:2}
Let $x$ be a real irrational algebraic number and $g\ge2$ a positive integer.
Then $x$ is normal in base~$g$.
\end{conjecture}

As shown by Borel~\cite{waldschmidt-borel09}, almost all numbers, in the sense of Lebesgue measure, are normal.
Examples of computable normal numbers have been constructed by Sierpinski, Lebesgue, Becher and Figueira; see~\cite{waldschmidt-BF02}.
However, the known algorithms to compute such examples are fairly complicated; indeed, ``ridiculously exponential'', according to~\cite{waldschmidt-BF02}.

An example of a $2$-normal number is the \emph{binary Champernowne number}, obtained by concatenation of the sequence of integers
\begin{displaymath}
0.\,1\,10\,11\,100\,101\,110\,111\,1000\,1001\,1010\,1011\,1100\,1101\,1110\,\ldots;
\end{displaymath}
see \cite{waldschmidt-champernowne33,waldschmidt-pillai39,waldschmidt-BC01}.
A closed formula for this number is
\begin{displaymath}
\sum_{k\ge1}k2^{-c_k},
\hspace{.15in}
c_k=k+\sum_{j=1}^k\lfloor\log_2j\rfloor;
\end{displaymath}
see \cite[page~183]{waldschmidt-BC01}.

Another example is given by Korobov, Stoneham, and others: if $a$ and $g$ are coprime integers greater than~$1$, then
\begin{displaymath}
\sum_{n\ge0}a^{-n}g^{-a^n}
\end{displaymath}
is normal in base~$g$; see~\cite{waldschmidt-BBCP04}.

A further example, due to Copeland and Erd\H{o}s in 1946, of a normal number in base $10$ is
\begin{displaymath}
0.23571113171923\ldots,
\end{displaymath}
obtained by concatenation of the sequence of prime numbers.
  
As pointed out to the author by Tanguy Rivoal, a definition of what it means for a number to have a
\emph{random sequence of digits} is given in~\cite{waldschmidt-calude02}, the first concrete example being \emph{Chaitin's omega number} which is the \emph{halting probability of a universal
self-delimiting computer with null-free data}.
This number, which is normal and transcendental, appears also to be a good candidate for being a
non-period in the sense of Kontsevich and Zagier~\cite{waldschmidt-KZ01}.

%%%%%%%%%%%%%%%%%%%%
%
% SUBSECTION 1.2
%
%%%%%%%%%%%%%%%%%%%%
  
\subsection{BBP Numbers}\label{SS:BBP}

An interesting approach towards Conjecture \ref{Conj:2} is provided by \emph{Hypothesis A} of Bailey and Crandall~\cite{waldschmidt-BC01}, who relate the question whether numbers like $\pi$, $\log2$ and other constants are normal to the following hypothesis involving the behaviour of the orbits of a discrete dynamical system.

\begin{hypothesisa}
Let
\begin{displaymath}
\theta:=\sum_{n\ge1}\frac{p(n)}{q(n)}\,g^{-n},
\end{displaymath}
where $g\ge2$ is an integer, $R=p/q\in\bfQ(X)$ a rational function with $q(n)\ne0$ for $n\ge1$ and
$\deg p<\deg q$.
Set $y_0=0$ and
\begin{displaymath}
y_{n+1}=gy_n+\frac{p(n)}{q(n)}\mod 1.
\end{displaymath}
Then the sequence $(y_n)_{n\ge1}$ either has finitely many limit points or is uniformly distributed modulo~$1$.
\end{hypothesisa}

A connection with special values of $G$-functions has been noted by
Lagarias~\cite{waldschmidt-lagarias01}.   
In his paper, Lagarias defines \emph{BBP numbers}, with reference to the paper
\cite{waldschmidt-BBP97} by Bailey, Borwein and Plouffe, as numbers of the form
\begin{displaymath}
\sum_{n\ge1}\frac{p(n)}{q(n)}\,g^{-n},
\end{displaymath}
where $g\ge2$ is an integer, $p$ and $q$ relatively prime polynomials in $\bfZ[X]$ with $q(n)\ne0$ for $n\ge1$.

Here are a few examples from~\cite{waldschmidt-BB01}.
Since
\begin{displaymath}
\sum_{n\ge1}\frac{1}{n}\,x^n=-\log(1-x)
\hspace{.15in}
\text{and}
\hspace{.15in}
\sum_{n\ge1}\frac{1}{2n-1}\,x^{2n-1}=\frac{1}{2}\log\frac{1+x}{1-x},
\end{displaymath}
it follows that $\log2$ is a BBP number in base $2$ as well as in base~$9$, since
\begin{equation}
\log2
=\sum_{n\ge1}\frac{1}{n}\,2^{-n}
=\sum_{n\ge1} \frac{6}{2n-1}\,3^{-2n}.
\label{E:log2}
\end{equation}
Also $\log3$ is a BBP number in base~$4$, since
\begin{displaymath}
\log3
=2\log2+\log\frac{3}{4}
=\sum_{n\ge0}\frac{1}{2^{2n}}\,\frac{1}{2n+1};
\end{displaymath}
and $\pi$ is a BBP number in base~$16$, since
\begin{equation}
\pi=\sum_{n\ge0}\left(\frac{4}{8n+1}-\frac{2}{8n+4}-\frac{1}{8n+5}-\frac{1}{8n+6}\right)2^{-4n}.
\label{E:pi}
\end{equation}
Further examples are $\pi^2$ in base $64$ and in base~$81$, $(\log 2)^2$ in base~$64$, and
$\zeta(3)$ in base $2^{12}=4096$.

%%%%%%%%%%%%%%%%%%%%
%
% SUBSECTION 1.3
%
%%%%%%%%%%%%%%%%%%%%

\subsection{Number of $1$'s in the Binary Expansion of a Real Number}

Denote by $B(x,n)$ the number of $1$'s among the first $n$ binary digits of an irrational real
number~$x$.

If $x,y$ and $x+y$ are positive and irrational, then for all sufficiently large~$n$,
\begin{displaymath}
B(x+y,n)\le B(x,n)+B(y,n)+1.
\end{displaymath}
If $x,y$ and $xy$ are positive and irrational, then for all sufficiently large~$n$,
\begin{displaymath}
B(xy,n)\le B(x,n)B(y,n)+\log_2\lfloor x+y+1\rfloor.
\end{displaymath}
If $x$ is positive and irrational, then for each integer $A>0$, the bound
\begin{displaymath}
B(x,n)B\left(\frac{A}{x},n\right)\ge n-1-\left\lfloor\log_2\left(x+\frac{A}{x}+1\right)\right\rfloor
\end{displaymath}
holds for all sufficiently large~$n$; see \cite{waldschmidt-BBCP04,waldschmidt-rivoal-a}.

A consequence is that if $a$ and $b$ are two integers, both at least~$2$, then none of the powers of the transcendental number
\begin{displaymath}
\xi=\sum_{n\ge1}a^{-b^n}
\end{displaymath}
is simply normal in base~$2$.
Also the lower bound
\begin{displaymath}
B(\sqrt{2},n)\ge n^{1/2}+O(1)
\end{displaymath}
can be deduced~\cite{waldschmidt-rivoal-a}.
In \cite[Theorem~7.1]{waldschmidt-BBCP04}, Bailey, Borwein, Crandall and Pomerance have obtained a similar lower bound valid for all real irrational algebraic numbers.

\begin{theorem}[Bailey, Borwein, Crandall, and Pomerance, 2004]\label{T:BBCP}
Let $x$ be a real algebraic number of degree at least~$2$.
Then there is a positive number~$C$, which depends only on~$x$, such that the number of $1$'s among the first $N$ digits in the binary expansion of $x$ is at least $CN^{1/d}$.
\end{theorem}

Further results related to Theorem~\ref{T:BBCP} are given by Rivoal~\cite{waldschmidt-rivoal-a}, Bugeaud~\cite{waldschmidt-bugeaud-a}, and Bugeaud and Evertse~\cite{waldschmidt-BE-a}.

As pointed out by Bailey, Borwein, Crandall and Pomerance, it follows from Theorem~\ref{T:BBCP} that for each $d\ge2$, the number
\begin{displaymath}
\sum_{n\ge0}2^{-d^n}
\end{displaymath}
is transcendental.
The transcendence of the number
\begin{equation}
\sum_{n\ge0}2^{-2^n}
\label{E:Kempner}
\end{equation}
goes back to Kempner in 1916; see \cite[Section~13.10]{waldschmidt-AS03}.
Other proofs are available, which rest either on Mahler's method, to be discussed in
Section~\ref{SS:Mahler}, or on the approximation theorem of Thue-Siegel-Roth-Ridout, to be discussed in Section~\ref{S:DA}; see also \cite[Section~1.6]{waldschmidt-schneider57}, and
\cite{waldschmidt-roth55,waldschmidt-roth60,waldschmidt-adamczewski04}.

Another consequence  of Theorem~\ref{T:BBCP} is the transcendence of the number
\begin{equation}
\sum_{n\ge0}g^{-F_n},
\label{E:puissanceFibonacci}
\end{equation}
for any integer $g\ge2$, where $(F_n)_{n\ge0}$ is the Fibonacci sequence, with $F_0=0$, $F_1=1$ and $F_{n+1}=F_n+F_{n-1}$ for $n\ge1$.
Again, the transcendence of this number also follows from Mahler's method \cite{waldschmidt-BBCP04} as well as from the theorem of Thue-Siegel-Roth-Ridout
\cite{waldschmidt-roth55,waldschmidt-roth60,waldschmidt-adamczewski04}.

%%%%%%%%%%%%%%%%%%%%
%
% SECTION 2
%
%%%%%%%%%%%%%%%%%%%%

\section{Words and Automata}\label{S:WordsAutomata}

%%%%%%%%%%%%%%%%%%%%
%
% SUBSECTION 2.1
%
%%%%%%%%%%%%%%%%%%%%

\subsection{Words}\label{SS:CoW}

We recall some basic facts from language theory; see for instance
\cite{waldschmidt-AS03,waldschmidt-lothaire02}.

We consider an alphabet $A$ with $g$ letters.
The free monoid $A^\ast$ on $A$ is the set of \emph{finite words} $a_1\ldots a_n$ where $n\ge0$ and $a_i\in A$ for $1\le i\le n$.
The law on $A^\ast$ is called \emph{concatenation}.

The number of letters of a finite word is its \emph{length}: the length of $a_1\ldots a_n$ is~$n$.

The number of words of length $n$ is $g^n$ for $n\ge0$.
The single word of length $0$ is the empty word $e$ with no letter.
It is the neutral element for the concatenation.

We shall consider \emph{infinite words} $w=a_1a_2a_3\ldots.$
A \emph{factor of length $m$} of such a word $w$ is a word of the form $a_ka_{k+1}\ldots a_{k+m-1}$ for some $k\ge1$.

The \emph{complexity} of an infinite word $w$ is the function $p(m)$ which counts, for each $m\ge1$, the number of distinct factors of $w$ of length~$m$.
Hence for an alphabet $A$ with $g$ elements we have $1\le p(m)\le g^m$, and the function
$m\mapsto p(m)$ is non-decreasing.
Conjecture \ref{Conj:1} is equivalent to the assertion that the complexity of the sequence of digits in base $g$ of an irrational algebraic number should be $p(m)=g^m$.

An infinite word is periodic if and only if its complexity is bounded.
If the complexity $p(m)$ of a word satisfies $p(m+1)=p(m)$ for one value of~$m$, then $p(m+k)=p(m)$ for all $k\ge0$, hence the word is periodic.
It follows that the complexity of a non-periodic word satisfies $p(m)\ge m+1$.
Following Morse and Hedlund, a word of minimal complexity $p(m)=m+1$ is called a \emph{Sturmian word}.
Sturmian words are those which encode with two letters the orbits of square billiard starting with an irrational angle; see \cite{waldschmidt-AMST94a,waldschmidt-AMST94b,waldschmidt-ST92}.
It is easy to check that on the  alphabet $\{a,b\}$, a Sturmian word $w$ is characterized by the property that for each $m\ge1$, there is exactly one factor $v$ of  $w$ of length $m$ for which both $va$ and
$vb$ are factors of $w$ of length $m+1$.

Let $A$ and $B$ be two finite sets.
A map from $A$ to $B^\ast$ can be uniquely extended to a homomorphism between the free monoids $A^\ast$ and~$B^\ast$.
We call such a homomorphism a \emph{morphism from $A$ to~$B$}.
The morphism is \emph{uniform} if all words in the image of $A$ have the same length.

Let $\varphi$ be a morphism from $A$ into itself.
Assume that there exists a letter $a$ for which $\varphi(a)=au$, where $u$ is a non-empty word satisfying  $\varphi^k(u)\ne e$ for every $k\ge0$.
Then the sequence of finite words $(\varphi^k(a))_{k\ge1}$ converges in~$A^{\bfN}$, endowed with the product topology of the discrete topology on each copy of~$A$, to an infinite word
$w=au\varphi(u)\varphi^2(u)\varphi^3(u)\ldots.$
This infinite word is clearly a fixed point for $\varphi$ and we say that \emph{$w$ is generated by the morphism $\varphi$}.

If, moreover, every letter occurring in $w$ occurs at least twice, then we say that \emph{$w$ is generated by a recurrent morphism}.

If the alphabet $A$ has two letters, then we say that $w$ is generated by a \emph{binary morphism}.

More generally, an infinite sequence $w$ in $A^{\bfN}$ is said to be \emph{morphic} (respectively
\emph{uniformly morphic}) if there exist a sequence $u$ generated by a morphism (respectively a uniform morphism) defined over an alphabet $B$ and a morphism (respectively a uniform morphism)
$\varphi$ from $B$ to $A$ such that $w=\varphi(u)$.

%%%%%%%%%%%%%%%%%%%%
%
% SUBSECTION 2.2
%
%%%%%%%%%%%%%%%%%%%%

\subsection{Finite Automata and Automatic Sequences}\label{SS:Automata}

A formal definition of a finite automaton is given, for instance, in \cite[Section~4.1]{waldschmidt-AS03}, \cite[Section 1.3.2]{waldschmidt-lothaire02}, \cite[Section~3.3]{waldschmidt-allouche05b} and
\cite[Section~3]{waldschmidt-AB07}.
We do not give the exact definition but we propose a number of examples in Section~\ref{SS:Examples}.
It suffices to say that a finite automaton consists of the following elements:
\begin{itemize}
\item the \emph{input alphabet}, which is usually the set of $g\ge2$ digits $\{0,1,\ldots,g-1\}$;
\item the \emph{set of states}~$\QQQ$, usually a finite set of $2$ or more elements, with one element, called the \emph{initial state} and denoted by~$i$, singled out; the elements of $\QQQ$ will be denoted by letters $\{i,a,b,\ldots\}$;
\item the \emph{transition map} $\QQQ\times\{0,1,\ldots,g-1\}\to\QQQ:(a,n)\mapsto n[a]$ which associates to every state $a$ a new state $n[a]$ depending on the current input~$n$;
\item the \emph{output alphabet}~$\AAA$, together with the \emph{output map} $\QQQ\to\AAA$.
\end{itemize}

We extend the transition map to a map $(a,n)\mapsto n[a]$ from $\QQQ\times\bfN\to\QQQ$ as follows: let $(a,n)\in\QQQ\times\bfN$, replace $n$ by its $g$-ary expansion $n=e_ke_{k-1}\ldots e_1e_0$ and define inductively $n[a]=e_ke_{k-1}\ldots e_1[b]$, where $b=e_0[a]$ is the image of $(a,e_0)$ under the transition map.
Hence an automaton produces an infinite sequence
\begin{displaymath}
0[i],1[i],2[i],3[i],\ldots,n[i],\ldots
\end{displaymath}
of elements of~$\QQQ$.
We take the images of the elements of this sequence under the transition map, and this gives rise to the \emph{output sequence}.
Such a sequence is called \emph{$g$-automatic}.
In other words, for an integer $g\ge2$, an infinite sequence $(a_n)_{n\ge0}$ of elements of
$\{0,1,\ldots,g-1\}$ is said to be \emph{$g$-automatic} if $a_n$ is a finite-state function of the representation of $n$ in base~$g$.

For instance let us explain why the characteristic function of the sequence $2^0=1$, $2^1=2$, $2^2=4$, $\ldots$ of powers of $2$ is $2$-automatic.
Take $g=2$, $\QQQ=\{i,a,b\}$, and define the transition map by
\begin{displaymath}
\begin{array}{ccc}
0[i]=i,
&0[a]=a,
&0[b]=b,\\
1[i]=a,
&1[a]=b,
&1[b]=b.
\end{array}
\end{displaymath}
A convenient notation is given by the following diagram.
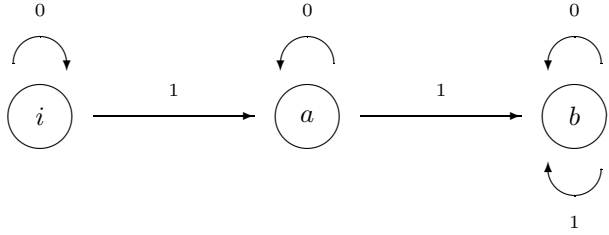
\begin{center}
\begin{picture}(300,100)(0,0)
\put(50,50){\circle{25}}
\put(50,50){\makebox(0,0){$i$}}
\put(150,50){\circle{25}}
\put(150,50){\makebox(0,0){$a$}}
\put(250,50){\circle{25}}
\put(250,50){\makebox(0,0){$b$}}
\put(70,50){\vector(1,0){60}}
\put(100,60){\makebox(0,0){$\scriptstyle1$}}
\put(170,50){\vector(1,0){60}}
\put(200,60){\makebox(0,0){$\scriptstyle1$}}
\put(50,70){\oval(20,20)[t]}
\put(60,67){\vector(0,-1){0}}
\put(50,90){\makebox(0,0){$\scriptstyle0$}}
\put(150,70){\oval(20,20)[t]}
\put(140,67){\vector(0,-1){0}}
\put(150,90){\makebox(0,0){$\scriptstyle0$}}
\put(250,70){\oval(20,20)[t]}
\put(240,67){\vector(0,-1){0}}
\put(250,90){\makebox(0,0){$\scriptstyle0$}}
\put(250,30){\oval(20,20)[b]}
\put(240,33){\vector(0,1){0}}
\put(250,10){\makebox(0,0){$\scriptstyle1$}}
\end{picture}
\end{center}
One easily checks, on writing the sequence of positive integers in base~$2$, that
\begin{align}
&
0[i]=i,
\hspace{.15in}
1[i]=a,
\hspace{.15in}
10[i]=a,
\hspace{.15in}
11[i]=b,
\nonumber
\\
&
100[i]=a,
\hspace{.15in}
101[i]=b,
\hspace{.15in}
110[i]=b,
\hspace{.15in}
\ldots.
\nonumber
\end{align}
Next define $f(i)=0$, $f(a)=1$ and $f(b)=0$. 
The output sequence
\begin{displaymath}
a_0a_1a_2\ldots=01101000100000001000\ldots
\end{displaymath}
is given by
\begin{displaymath}
a_n=\left\{\begin{array}{ll}
1&\text{if $n$ is a power of $2$},\\
0&\text{otherwise}.
\end{array}
\right.
\end{displaymath}

According to Cobham~\cite{waldschmidt-cobham72}, automatic sequences have complexity
$p(m)=O(m)$; see also \cite[Section~10.3]{waldschmidt-AS03}.
Also, automatic sequences are the same as \emph{uniform morphic sequences}; see for instance
\cite[Section~6.3]{waldschmidt-AS03} and \cite[Theorem~4.1]{waldschmidt-AC06}.

Automatic sequences are between periodicity and chaos.
They occur in connection with harmonic analysis, ergodic theory, fractals, Feigenbaum trees,
quasi-crystals, transition phases in statistical mechanics, and others; see
\cite{waldschmidt-AM88,waldschmidt-cerf04,waldschmidt-allouche05b} and
\cite[Chapter~17]{waldschmidt-AS03}.

%%%%%%%%%%%%%%%%%%%%
%
% SUBSECTION 2.3
%
%%%%%%%%%%%%%%%%%%%%

\subsection{Examples}\label{SS:Examples}

%%%%%%%%%%%%%%%%%%%%
%
% SUBSUBSECTION 2.3.1
%
%%%%%%%%%%%%%%%%%%%%

\subsubsection{The Fibonacci Word}\label{SS:FW}

Consider the alphabet $A=\{a,b\}$.
Start with $f_1=b$, $f_2=a$ and define $f_n=f_{n-1}f_{n-2}$.
Then
\begin{displaymath}
f_3=ab,
\hspace{.1in}
f_4=aba,
\hspace{.1in}
f_5=abaab,
\hspace{.1in}
f_6=abaababa,
\hspace{.1in}
f_7=abaababaabaab,
\hspace{.1in}
\ldots.
\end{displaymath}
There is a unique word
\begin{displaymath}
w=abaababaabaababaababaabaababaabaab\ldots
\end{displaymath}
in which $f_n$ is the prefix of length~$F_n$, the Fibonacci number of index~$n$, for $n\ge2$.
This is the \emph{Fibonacci word}; it  is generated by a binary recurrent morphism
\cite[Section~7.1]{waldschmidt-AS03} and is the fixed point of the morphism $a\mapsto ab$,
$b\mapsto a$; under this morphism, the image of $f_n$ is $f_{n+1}$.

Let us check that the Fibonacci word is not periodic.
Indeed, the word $f_n$ has length~$F_n$, and consists of $F_{n-1}$ letters $a$ and $F_{n-2}$
letters~$b$. 
Hence the proportion of $a$ in the Fibonacci word $w$ is $1/\Phi$, where $\Phi$ is  the golden number
\begin{displaymath}
\Phi=\frac{1+\sqrt{5}}{2}
\end{displaymath}
which is an irrational number.
On the other hand, for a periodic word the proportion of each letter is a rational number.

\begin{remark}
The proportion of $b$ in $w$ is $1/\Phi^2$, with $1/\Phi+1/\Phi^2=1$, as expected!
\end{remark}

\begin{proposition}\label{P:FibonacciSturmian}
The Fibonacci word is Sturmian.
\end{proposition}

We write the factors  of length $1,2,3,4,5,\ldots$ of the Fibonacci word in columns as shown.
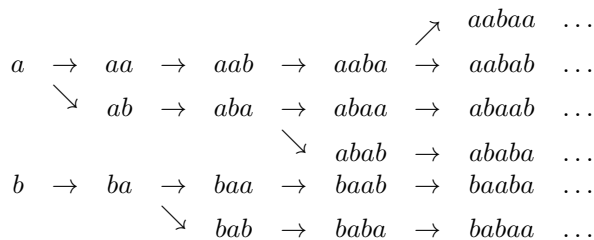
\begin{displaymath}
\begin{array}{cccccccccc}
&&&&&&&&aabaa&\ldots\\[-.1in]
&&&&&&&\nearrow\\
a&\to&aa&\to&aab&\to&aaba&\to&aabab&\ldots\\
&\searrow\\[-.1in]
&&ab&\to&aba&\to&abaa&\to&abaab&\ldots\\
&&&&&\searrow\\[-.1in]
&&&&&&abab&\to&ababa&\ldots\\
b&\to&ba&\to&baa&\to&baab&\to&baaba&\ldots\\
&&&\searrow\\[-.1in]
&&&&bab&\to&baba&\to&babaa&\ldots
\end{array}
\end{displaymath}
The $k$-th column contains $k+1$ elements: one of them has two right extensions to a factor of length $k+1$, and all remaining $k$ factors of length $k$ have a single extension.
This is easily seen to be a characterization of Sturmian words.

Concerning the sequence
\begin{equation}
(v_n)_{n\ge0}=(0,1,0,0,1,0,1,0,0,1,0,0,1,0,\ldots),
\label{E:Danilov}
\end{equation}
derived from the Fibonacci word on the alphabet $\{0,1\}$, a result of Danilov
\cite{waldschmidt-danilov72} in 1972 states that for all integers $g\ge2$, the number
\begin{displaymath}
\sum_{n\ge0}v_ng^{-n}
\end{displaymath}
is transcendental; see also \cite[Theorem~4.2]{waldschmidt-allouche05a}.

\begin{proposition}\label{P:FibonacciNotAutomatic}
The Fibonacci word is not automatic.
\end{proposition}
 
Proposition~\ref{P:FibonacciNotAutomatic} follows from a result of Cobham
\cite{waldschmidt-cobham72} which shows that the frequency of a letter in an automatic word, if it exists, is a rational number. 

The origin of the Fibonacci sequence is a model of growth of a population of rabbits.
Denote a pair of young rabbits by $Y$ and a pair of adults by~$A$.
From one year to the next, the young pair become adult, which we write as $Y\to A$, while the adult pair stay alive and produce a young pair, which we write as $A\to AY$.
This gives rise to the dynamical system
\begin{displaymath}
Y
\to A
\to AY
\to AYA
\to AYA\,AY
\to AYA\,AYAYA
\to\ldots,
\end{displaymath}
and the sequence $(R_1,R_2,R_3,\ldots)$ of Fibonacci rabbits
\begin{displaymath}
A,Y,A,A,Y,A,Y,A,A,Y,A,A,Y,A,\ldots.
\end{displaymath}
Replacing $Y$ by $1$ and $A$ by $0$ produces the sequence \eqref{E:Danilov} considered by Danilov.

Any integer $n\ge2$ has a unique representation as the sum of distinct Fibonacci numbers $F_m$, $m\ge2$, with the property that no Fibonacci numbers with consecutive indices occur in the sum.
This representation yields the following algorithm to decide whether $R_n$ is $A$ or~$Y$.
If the smallest index in the decomposition of $n$ is even, then $R_n=A$.
If it is odd, then $R_n=Y$.
For instance, for $51=F_9+F_7+F_4+F_2$, the smallest index, namely~$2$, is even, so $R_{51}=A$.

Recall that $\Phi=(1+\sqrt{5})/2$ denotes the golden number.
The sequence of indices $n$ for which $R_n=A$ is
\begin{displaymath}
\lfloor\Phi\rfloor=1,
\hspace{.1in}
\lfloor2\Phi\rfloor=3,
\hspace{.1in}
\lfloor3\Phi\rfloor=4,
\hspace{.1in}
\lfloor4\Phi\rfloor=6,
\hspace{.1in}
\ldots,
\end{displaymath}
while the sequence of indices $n$ for which $R_n=Y$ is
\begin{displaymath}
\lfloor\Phi^2\rfloor=2,
\hspace{.1in}
\lfloor2\Phi^2\rfloor=5,
\hspace{.1in}
\lfloor3\Phi^2\rfloor=7,
\hspace{.1in}
\lfloor4\Phi^2\rfloor=8,
\hspace{.1in}
\ldots.
\end{displaymath}
For instance, $32\Phi=51.77\ldots$, so $\lfloor32\Phi\rfloor=51$ and $R_{51}=A$.

This sequence $(R_n)_{n\ge1}$ of rabbits is an example of a \emph{Beatty sequence}.

%%%%%%%%%%%%%%%%%%%%
%
% SUBSUBSECTION 2.3.2
%
%%%%%%%%%%%%%%%%%%%%

\subsubsection{The Prouhet-Thue-Morse word $abbabaabbaababbab\ldots\,$}

The finite automaton
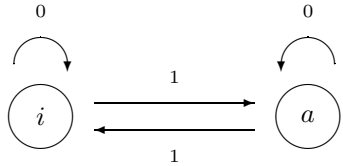
\begin{center}
\begin{picture}(200,75)(0,0)
\put(50,25){\circle{25}}
\put(50,25){\makebox(0,0){$i$}}
\put(150,25){\circle{25}}
\put(150,25){\makebox(0,0){$a$}}
\put(70,30){\vector(1,0){60}}
\put(100,40){\makebox(0,0){$\scriptstyle1$}}
\put(130,20){\vector(-1,0){60}}
\put(100,10){\makebox(0,0){$\scriptstyle1$}}
\put(50,45){\oval(20,20)[t]}
\put(60,42){\vector(0,-1){0}}
\put(50,65){\makebox(0,0){$\scriptstyle0$}}
\put(150,45){\oval(20,20)[t]}
\put(140,42){\vector(0,-1){0}}
\put(150,65){\makebox(0,0){$\scriptstyle0$}}
\end{picture}
\end{center}
with $f(i)=0$ and $f(a)=1$ produces the sequence $a_0a_1a_2\ldots,$ where $a_n=f(n[i])$.
For instance, with $n=9$, since $9$ is $1001$ in binary notation,
\begin{displaymath}
1001[i]
=100[a]
=10[a]
=1[a]
=i,
\end{displaymath}
and $f(i)=0$, we have $a_9=0$.

This is the \emph{Prouhet-Thue-Morse sequence}
\begin{displaymath}
01101001100101101\ldots,
\end{displaymath}
where the $(n+1)$-th term $a_n$ is $1$ if the number of $1$'s, which is the same as the sum of the binary digits, in the binary expansion of $n$ is odd, and is $0$ otherwise; see
\cite[Section~1.6]{waldschmidt-AS03}.

If, in the Prouhet-Thue-Morse sequence, we replace $0$ by $a$ and $1$ by~$b$, we obtain the
\emph{Prouhet-Thue-Morse word} on the alphabet $\{a,b\}$, and this starts with
\begin{displaymath}
w=abbabaabbaababbab\ldots.
\end{displaymath}
This word is generated by a binary recurrent morphism; see \cite[Section~6.2]{waldschmidt-AS03}.
It is the fixed point of the morphism $a\mapsto ab$, $b\mapsto ba$.

An interesting observation, due to Thue in 1906, is that if $w$ is a finite word and $a$ a letter for which $wwa$ is a factor of the Prouhet-Thue-Morse word, then $a$ is not the first letter of~$w$.
Therefore, in the Prouhet-Thue-Morse sequence, no three consecutive identical blocks such as $000$ or $111$ or $010101$ or  $101010$ or $001001001$ can occur.

%%%%%%%%%%%%%%%%%%%%
%
% SUBSUBSECTION 2.3.3
%
%%%%%%%%%%%%%%%%%%%%

\subsubsection{The Baum-Sweet Sequence}

For $n\ge0$, let $a_n=1$ if the binary expansion of $n$ contains no block of consecutive $0$'s of odd length, and $a_n=0$ otherwise.
The \emph{Baum-Sweet sequence} $(a_n)_{n\ge0}$ starts with
\begin{displaymath}
110110010100100110010\ldots.
\end{displaymath}
This sequence is automatic, associated with  the automaton 
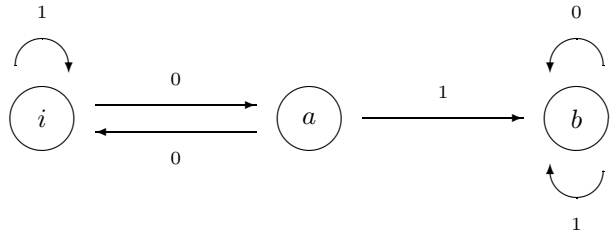
\begin{center}
\begin{picture}(300,100)(0,0)
\put(50,50){\circle{25}}
\put(50,50){\makebox(0,0){$i$}}
\put(150,50){\circle{25}}
\put(150,50){\makebox(0,0){$a$}}
\put(250,50){\circle{25}}
\put(250,50){\makebox(0,0){$b$}}
\put(70,55){\vector(1,0){60}}
\put(100,65){\makebox(0,0){$\scriptstyle0$}}
\put(130,45){\vector(-1,0){60}}
\put(100,35){\makebox(0,0){$\scriptstyle0$}}
\put(170,50){\vector(1,0){60}}
\put(200,60){\makebox(0,0){$\scriptstyle1$}}
\put(50,70){\oval(20,20)[t]}
\put(60,67){\vector(0,-1){0}}
\put(50,90){\makebox(0,0){$\scriptstyle1$}}
\put(250,70){\oval(20,20)[t]}
\put(240,67){\vector(0,-1){0}}
\put(250,90){\makebox(0,0){$\scriptstyle0$}}
\put(250,30){\oval(20,20)[b]}
\put(240,33){\vector(0,1){0}}
\put(250,10){\makebox(0,0){$\scriptstyle1$}}
\end{picture}
\end{center}
with $f(i)=1$, $f(a)=0$ and $f(b)=0$; see \cite[Example 5.1.7]{waldschmidt-AS03}.

%%%%%%%%%%%%%%%%%%%%
%
% SUBSUBSECTION 2.3.4
%
%%%%%%%%%%%%%%%%%%%%

\subsubsection{Powers of $2$}

As we have seen, the binary number
\begin{displaymath}
\xi
:=\sum_{n\ge0}2^{-2^n}
=0.1101000100000001000\ldots
=0.a_1a_2a_3\ldots
\end{displaymath}
is $2$-automatic.
The associated infinite word
\begin{displaymath}
\bfv
=v_1v_2v_3\ldots
=bbabaaabaaaaaaabaaa\ldots,
\end{displaymath}
where
\begin{displaymath}
v_n=\left\{\begin{array}{ll}
b&\text{if $n$ is a power of $2$},\\
a&\text{otherwise},
\end{array}
\right.
\end{displaymath}
has complexity $p(m)$ bounded by~$2m$; the initial values are given below.
\begin{displaymath}
\begin{array}{ccccccccc}
m&&1&2&3&4&5&6&\ldots\\
p(m)&&2&4&6&7&9&11&\ldots
\end{array}
\end{displaymath}
%
%

%%%%%%%%%%%%%%%%%%%%
%
% SUBSUBSECTION 2.3.5
%
%%%%%%%%%%%%%%%%%%%%

\subsubsection{The Rudin-Shapiro Word}

For $n\ge0$, let $r_n\in\{a,b\}$ satisfy $r_n=a$ (respectively $r_n=b$) if the number of occurrences of 
the pattern $11$ in the binary representation of $n$ is even (respectively odd).
This produces the \emph{Rudin-Shapiro word}
\begin{displaymath}
aaabaabaaaabbbab\ldots.
\end{displaymath}

Let $\sigma$ be the morphism defined from the monoid $B^\ast$ on the alphabet $B=\{1,2,3,4\}$ into $B^\ast$ by $\sigma(1)=12$, $\sigma(2)=13$, $\sigma(3)=42$ and $\sigma(4)=43$.
Let
\begin{displaymath}
\bfu=121312421213\ldots
\end{displaymath}
be the fixed point of $\sigma$ beginning with $1$ and let $\varphi$ be the morphism defined from
$B^\ast$ to $\{a,b\}^\ast$ by $\varphi(1)=aa$, $\varphi(2)= ab$, $\varphi(3)=ba$ and $\varphi(4)=bb$.
Then the \emph{Rudin-Shapiro word} is $\varphi(\bfu)$, hence it is morphic.

%%%%%%%%%%%%%%%%%%%%
%
% SUBSUBSECTION 2.3.6
%
%%%%%%%%%%%%%%%%%%%%

\subsubsection{Paper Folding}

Folding  a strip of paper always in the same direction, and then opening it up, yields a sequence of folds which can be encoded using the digits $0$ and~$1$.
The resulting sequence $(u_n)_{n\ge0}$, given by
\begin{displaymath}
1101100111001001\ldots,
\end{displaymath}
satisfies
\begin{displaymath}
u_{4n}=1,
\hspace{.15in}
u_{4n+2}=0,
\hspace{.15in}
u_{2n+1}=u_n,
\end{displaymath}
and is produced by the automaton
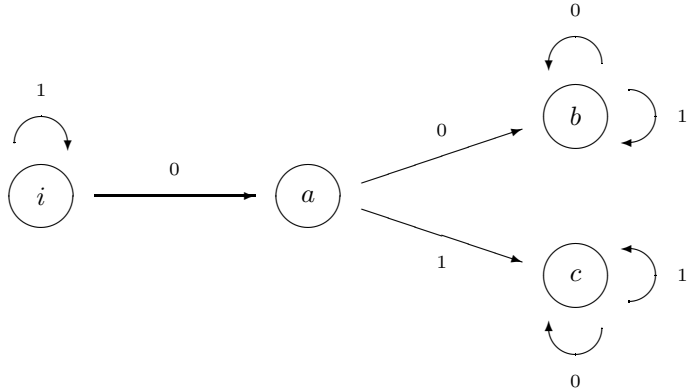
\begin{center}
\begin{picture}(330,160)(0,0)
\put(50,80){\circle{25}}
\put(50,80){\makebox(0,0){$i$}}
\put(150,80){\circle{25}}
\put(150,80){\makebox(0,0){$a$}}
\put(250,110){\circle{25}}
\put(250,110){\makebox(0,0){$b$}}
\put(250,50){\circle{25}}
\put(250,50){\makebox(0,0){$c$}}
\put(70,80){\vector(1,0){60}}
\put(100,90){\makebox(0,0){$\scriptstyle0$}}
\put(170,85){\vector(3,1){60}}
\put(200,105){\makebox(0,0){$\scriptstyle0$}}
\put(170,75){\vector(3,-1){60}}
\put(200,55){\makebox(0,0){$\scriptstyle1$}}
\put(50,100){\oval(20,20)[t]}
\put(60,97){\vector(0,-1){0}}
\put(50,120){\makebox(0,0){$\scriptstyle1$}}
\put(250,130){\oval(20,20)[t]}
\put(240,127){\vector(0,-1){0}}
\put(250,150){\makebox(0,0){$\scriptstyle0$}}
\put(270,110){\oval(20,20)[r]}
\put(267,100){\vector(-1,0){0}}
\put(290,110){\makebox(0,0){$\scriptstyle1$}}
\put(250,30){\oval(20,20)[b]}
\put(240,33){\vector(0,1){0}}
\put(250,10){\makebox(0,0){$\scriptstyle0$}}
\put(270,50){\oval(20,20)[r]}
\put(267,60){\vector(-1,0){0}}
\put(290,50){\makebox(0,0){$\scriptstyle1$}}
\end{picture}
\end{center}
with $f(i)=f(a)=f(b)=1$ and $f(c)=0$.
 
An equivalent definition for this sequence is given as follows: the sequence $a_n=u_{n+1}$, $n\ge1$, is defined recursively by $a_n=1$ if $n$ is a power of~$2$, and
\begin{displaymath}
a_{2^k+a}=1-a_{2^k-a},
\hspace{.15in}
1\le a<2^k;
\end{displaymath}
see \cite[Example 5.1.6]{waldschmidt-AS03}.

For a connection between the paper folding sequence and the Prouhet-Thue-Morse sequence,
see~\cite{waldschmidt-bacher06}.

%%%%%%%%%%%%%%%%%%%%
%
% SUBSECTION 2.4
%
%%%%%%%%%%%%%%%%%%%%

\subsection{Complexity of the $g$-ary Expansion of an Algebraic Number}

The transcendence of a number whose sequence of digits is Sturmian has  been proved by Ferenczi and Mauduit \cite{waldschmidt-FM97} in 1997.
The point is that such sequences contain sequences of digits which bear similarities, and yields the existence of very sharp rational approximations which do not exist for algebraic numbers.

It follows from their work that the complexity of the $g$-ary expansion of every irrational algebraic number satisfies
\begin{displaymath}
\liminf_{m\to\infty}(p(m)-m)=+\infty;
\end{displaymath}
see~\cite{waldschmidt-allouche00}.
The main tool for the proof is a $p$-adic version of the Thue-Siegel-Roth theorem due to Ridout in 1957; see Theorem~\ref{T:Ridout} below as well as~\cite{waldschmidt-AB05b}.

Several papers have been devoted to the study of the complexity of the $g$-ary expansions of real algebraic numbers, in particular by Allouche and Zamboni in 1998, Risley and Zamboni in 2000, and Adamczewski and Cassaigne in 2003.
For a survey, see~\cite{waldschmidt-AB05b}. 
The main recent result is the following~\cite{waldschmidt-AB07}.

\begin{theorem}[Adamczewski and Bugeaud, 2007]\label{T:ABL}
The complexity $p(m)$ of the $g$-ary expansion of a real irrational algebraic number satisfies
\begin{displaymath}
\liminf_{m\to\infty}\frac{p(m)}{m}=+\infty.
\end{displaymath}
\end{theorem}

In 1968\footnote{The author is grateful to Boris Adamczewski for this reference.}, Cobham
\cite{waldschmidt-cobham68} claimed that automatic irrational numbers are transcendental.
This follows from Theorem~\ref{T:ABL}, since automatic numbers have a complexity $O(m)$.

\begin{corollary}[Conjecture of Cobham]\label{C:Cobham}
If the sequence of digits of an irrational real number $x$ is automatic, then $x$ is transcendental.
\end{corollary}

The main tool for the proof of Theorem~\ref{T:ABL} is a new, combinatorial transcendence criterion obtained by Adamczewski, Bugeaud and Luca \cite{waldschmidt-ABF04} as an application of Schmidt's subspace theorem; see Theorem~\ref{T:SST} and~\cite{waldschmidt-schmidt80}.

In 1979, Christol \cite{waldschmidt-christol79} proved that when $p$ is a prime and
\begin{displaymath}
f(X)=\sum_{k\ge1}u_kX^k\in\bfF_p[[X]]
\end{displaymath}
a power series with coefficients in the finite field~$\bfF_p$, then $f$ is algebraic if and only if the sequence of coefficients $(u_k)_{k\ge1}$ can be produced by a $p$-automaton.

According to Cobham, a sequence can be simultaneously $k$-automatic and $\ell$-automatic with $k$ and $\ell$ multiplicatively independent if and only if it is ultimately periodic.
Christol, Kamae, Mend\`es France and Rauzy  \cite{waldschmidt-CKMR80} have shown that a power series with coefficients in a finite set can be algebraic over two different finite fields if and only if it is rational.

In this context Theorem~\ref{T:ABL} implies the following statement.

\begin{corollary}\label{Corollaire}
Let $p$ be a prime number, $g\ge p$ an integer and $(u_k)_{k\ge1}$ a sequence of integers in the range $\{0,\ldots,p-1\}$.
The formal power series
\begin{displaymath}
\sum_{k\ge1}u_kX^k
\end{displaymath}
and the real number
\begin{displaymath}
\sum_{k\ge1}u_kg^{-k}
\end{displaymath}
are simultaneously algebraic over $\bfF_p(X)$ and over $\bfQ$ respectively if and only if they are rational.
\end{corollary}

As an example, taken from \cite[Section~6]{waldschmidt-allouche05a} and
\cite[Section~2.4]{waldschmidt-allouche05b}), consider the Prouhet-Thue-Morse sequence
$(a_n)_{n\ge0}$.
The series
\begin{displaymath}
F(X)=\sum_{n\ge0}a_nX^n
\end{displaymath}
is algebraic over $\bfF_2(X)$, as it is a root of $(1+X)^3F^2+(1+X)^2F+X=0$.
Since $F$ is not a rational function, Corollary~\ref{Corollaire} gives another proof of Mahler's transcendence result on the number
\begin{displaymath}
\sum_{n\ge0}a_ng^{-n}.
\end{displaymath}

The following result related to Theorem~\ref{T:ABL} has been obtained subsequently by Bugeaud and Evertse \cite{waldschmidt-BE-a} by means of a refinement of the so-called Cugiani-Mahler theorem.
The main tool is again a quantitative version of Schmidt's subspace
theorem~\cite{waldschmidt-schmidt80}.

\begin{theorem}[Bugeaud and Evertse, 2007]\label{T:BE}
Let $b\ge2$ be an integer and $\xi$ an irrational algebraic number with $0<\xi<1$.
Then for any real number $\eta< 1/11$, the complexity $p(m)$ of the $b$-ary expansion of $\xi$ satisfies
\begin{displaymath}
\limsup_{m\to+\infty}\frac{p(m)}{m(\log m)^{\eta}}=+\infty.
\end{displaymath}
\end{theorem}

Further developments of the method of \cite{waldschmidt-ABF04,waldschmidt-AB07} have recently been achieved by Adamczewski and Rampersad \cite{waldschmidt-AR-a} who have shown that the binary expansion of an algebraic number contains infinitely many occurrences of $7/3$-powers.
Here, for a positive real number~$w$, a \emph{$w$-th power} is a word of the form
$W^{\lfloor w\rfloor}W'$, where $W$ is a word of length~$\vert W\vert$, say, and $W'$ is the prefix of $W$ of length $\lceil(w-\lfloor w\rfloor)\vert W\vert\rceil$.
For instance with $w=7/3$ and $W=011010$ we have $\vert W\vert=6$, $\lfloor w\rfloor=2$ and
$\lceil(w-\lfloor w\rfloor)\vert W\vert\rceil=2$, so $011010\,011010\,01$ is a $7/3$-th power.

Adamczewski and Rampersad deduce from their result that the binary expansion of an algebraic number contains infinitely many occurrences of overlaps, where an overlap is a pattern of the form $xXxXx$, where $x$ is a letter and $X$ a word.

%%%%%%%%%%%%%%%%%%%%
%
% SECTION 3
%
%%%%%%%%%%%%%%%%%%%%

\section{Analytic Methods}\label{S:AnalyticMethods}

%%%%%%%%%%%%%%%%%%%%
%
% SUBSECTION 3.1
%
%%%%%%%%%%%%%%%%%%%%

\subsection{Mahler's Method}\label{SS:Mahler}

Mahler \cite{waldschmidt-mahler30} initiated a new transcendence method in 1929 for studying values of functions satisfying certain functional equations.
His work was not widely known and in 1969 he published \cite{waldschmidt-mahler69} which was the source of a revival of this method; see \cite{waldschmidt-LP77,waldschmidt-nishioka96}.

A first example \cite[Theorem 1.1.2]{waldschmidt-nishioka96} is the function
\begin{displaymath}
f(z)=\sum_{n\ge0}z^{-d^n},
\hspace{.15in}
d\ge2,
\end{displaymath}
which satisfies the functional equation$f(z^d)+z=f(z)$ for $\vert z\vert<1$.
Mahler proved that it takes transcendental values at algebraic points in the domain $0<\vert z\vert<1$.
A special case is \eqref{E:Kempner}.
     
Mahler also proved in 1929 that the so-called \emph{Prouhet-Thue-Morse-Mahler number in base
$g\ge2$}, given by
\begin{displaymath}
\xi_g=\sum_{n\ge0}\frac{a_n}{g^n},
\end{displaymath}
where$(a_n)_{n\ge0}$ is the Prouhet-Thue-Morse sequence, is transcendental; see
\cite{waldschmidt-nishioka96} and \cite[Section~13.4]{waldschmidt-AS03}.
The idea of proof is as follows; see \cite[Section~13.4]{waldschmidt-AS03} and
\cite[Example 1.3.1]{waldschmidt-nishioka96}, where the complete proof is given.
Consider the function
\begin{displaymath}
f(z)=\prod_{n\ge0}(1-z^{2^n})
\hspace{.15in}
\text{which satisfies}
\hspace{.15in}
f(z)=\sum_{n\ge0}(-1)^{a_n}z^n.
\end{displaymath}
For $a\in\{0,1\}$, we can write $(-1)^a=1-2a$.
Hence
\begin{displaymath}
f(z)
=\sum_{n\ge0}(1-2a_n)z^n
=\frac{1}{1-z}-2\sum_{n\ge0}a_nz^n.
\end{displaymath}
Using the functional equation $f(z)=(1-z)f(z^2)$, Mahler proves that $f(\alpha)$ is transcendental for all algebraic numbers $\alpha$ satisfying $0<\vert\alpha\vert<1$.

Another application \cite{waldschmidt-BBCP04} of Mahler's method is the transcendence of the number \eqref{E:puissanceFibonacci} whose digits in a given base $g$ are $1$ at the Fibonacci indices $1,2,3,5,8,\ldots,$ and $0$ elsewhere.

After Cobham \cite{waldschmidt-cobham72}, Loxton and van der Poorten also tried\footnote{See
\cite{waldschmidt-LP77} and \cite[Section~13.10]{waldschmidt-AS03}.} in 1982 and 1988 to use Mahler's method to prove Cobham's conjecture, now Corollary~\ref{C:Cobham} of Theorem~\ref{T:ABL} of Adamczewski and Bugeaud.
Becker \cite{waldschmidt-becker94} pointed out in 1994 that Mahler's method yields only a weaker result so far, that for any given non-eventually periodic automatic sequence $\bfu=(u_1,u_2,u_3\ldots)$, the real number
\begin{displaymath}
\sum_{k\ge1}u_kg^{-k}
\end{displaymath}
is transcendental, provided that the integer $g$ is sufficiently large in terms of~$\bfu$; see
also~\cite{waldschmidt-AB07}.
It is yet a challenge to extend Mahler's method in order to prove Cobham's conjecture.

There is a further very interesting development of Mahler's method which we only allude to here without entering the subject.
It is due to Denis and deals with transcendence problems in finite characteristic; see Pellarin's report \cite{waldschmidt-pellarin07} in the Bourbaki Seminar.

%%%%%%%%%%%%%%%%%%%%
%
% SUBSECTION 3.2
%
%%%%%%%%%%%%%%%%%%%%

\subsection{Nesterenko's Theorem and Consequences}\label{SS:Nesterenko}

The transcendence of the Liouville-Fredholm number \eqref{E:Liouville--Fredholm} below was studied only ten years ago by Bertrand in 1997, and Duverney, Nishioka as well as Nishioka and Shiokawa in 1998, as a consequence of the results due to Nesterenko's work in 1996 on the transcendence of values of theta series at rational points involving modular functions; see
\cite[Chapter~3, Section~1.3]{waldschmidt-NP01}.
It follows from these works that for algebraic $\alpha$ in the domain $0<\vert\alpha\vert<1$, the three numbers
\begin{displaymath}
\sum_{n\ge0}\alpha^{-n^2},
\hspace{.15in}
\sum_{n\ge0}n^2\alpha^{-n^2},
\hspace{.15in}
\sum_{n\ge0}n^4\alpha^{-n^2}
\end{displaymath}
are algebraically independent. 

A related example from \cite{waldschmidt-AB07} is the number
\begin{displaymath}
\eta=\sum_{k\ge1}u_k3^{-k}
\end{displaymath}
which is associated to the word
\begin{displaymath}
\bfu=012^112^212^312^412^512^612^7\ldots,
\end{displaymath}
where, for instance, $2^4$ denotes $2222$, generated by the non-recurrent morphism $0\mapsto012$, $1\mapsto12$, $2\mapsto2$.
Using
\begin{displaymath}
\eta=1-\sum_{n\ge1}3^{-n(n+1)/2},
\end{displaymath}
one deduces from Nesterenko's result that it is transcendental.
For this number the growth of complexity $p(m)$ is quadratic in~$m$, so Theorem~\ref{T:ABL} does not apply.

%%%%%%%%%%%%%%%%%%%%
%
% SECTION 4
%
%%%%%%%%%%%%%%%%%%%%

\section{Diophantine Approximation}\label{S:DA}

%%%%%%%%%%%%%%%%%%%%
%
% SUBSECTION 4.1
%
%%%%%%%%%%%%%%%%%%%%

\subsection{Liouville, Thue, Siegel, Roth, Ridout, Schmidt}\label{SS:DA}

Diophantine approximation theory yields information on the arithmetic nature of numbers of the form
\begin{equation}
\sum_{n\ge0}g^{-u_n}.
\label{E:serie}\end{equation}
For instance, if the sequence satisfies
\begin{displaymath}
u_{n+1}-u_n\to+\infty,
\hspace{.15in}
n\to+\infty,
\end{displaymath}
then the number given by \eqref{E:serie} is irrational.
This follows from the fact that  a real number is rational if and only if its $g$-ary expansion is ultimately periodic. It also follows from Diophantine analysis.
As a special case take $u_n=n^2$.
The irrationality of the Liouville-Fredholm number
\begin{equation}
\theta=\sum_{n\ge0}g^{-n^2}
\label{E:Liouville--Fredholm}\end{equation}  
follows.

Here is the very simple proof, which goes back to Liouville in 1851; see
\cite[Chapter~3, Section~1.3]{waldschmidt-NP01}.
Let $\epsilon>0$.
Let $N$ be a sufficiently large integer.
Set
\begin{displaymath}
q=a^{N^2}
\hspace{.15in}
\text{and}
\hspace{.15in}
p=\sum_{n=0}^Na^{N^2-n^2}.
\end{displaymath}
Then
\begin{displaymath}
0
<q\theta-p
=\sum_{n\ge N+1}\frac{1}{a^{n^2-N^2}}
=\sum_{k\ge1}\frac{1}{a^{2Nk+k^2}}
<\frac{\theta-1}{a^{2N}}.
\end{displaymath}
For $N$ sufficiently large the right hand side is less than $\epsilon$ and the irrationality of $\theta$ is proved.

That $\theta$ is not a quadratic irrationality follows from the work of Bailey and
Crandall~\cite{waldschmidt-BC01}.
See Section~\ref{SS:Nesterenko} for the transcendence of this number~$\theta$.

The first transcendence statement on numbers given by a series \eqref{E:serie}, assuming that the   sequence $(u_n)_{n\ge0}$ is increasing and grows sufficiently fast, goes back to Liouville in 1844.

\begin{theorem}[Liouville, 1844]
For any real algebraic number~$\alpha$, there exists a constant $c>0$ such that the set of $p/q\in\bfQ$ with $\vert\alpha-p/q\vert<q^{-c}$ is finite.
\end{theorem}

Liouville's theorem yields the transcendence of the value of a series like \eqref{E:serie}, provided that the increasing sequence $(u_n)_{n\ge0}$ satisfies
\begin{displaymath}
\limsup_{n\to\infty}\frac{u_{n+1}}{u_n}=+\infty.
\end{displaymath}
For instance, $u_n=n!$ satisfies this condition, so the number
\begin{displaymath}
\sum_{n\ge0}g^{-n!}
\end{displaymath}
is transcendental.

\begin{theorem}[Thue, Siegel, Roth \cite{waldschmidt-roth55,waldschmidt-roth60}, 1955]\label{T:Roth}
For any real algebraic number $\alpha$ and any $\epsilon>0$, the set of $p/q\in\bfQ$ with
$\vert\alpha-p/q\vert<q^{-2-\epsilon}$ is finite.
\end{theorem}
   
From Theorem~\ref{T:Roth}, one deduces the transcendence of the series \eqref{E:serie} under the weaker hypothesis
\begin{displaymath}
\limsup_{n\to\infty}\frac{u_{n+1}}{u_n}>2.
\end{displaymath}

The sequence $u_n=\lfloor\kappa^n\rfloor$ satisfies this condition as soon as $\kappa>2$.
For example, the transcendence of the number
\begin{displaymath}
\sum_{n\ge0}g^{-3^n}
\end{displaymath}
follows from Theorem~\ref{T:Roth}.

A stronger result follows from Theorem~\ref{T:Ridout} below due to Ridout (see, for instance,
\cite[Section~13.5]{waldschmidt-AS03}), using the fact that the denominators $g^{u_n}$ are powers
of~$g$.
The condition
\begin{displaymath}
\limsup_{n\to\infty}\frac{u_{n+1}}{u_n}>1
\end{displaymath}
is sufficient to imply the transcendence of the sum of the series \eqref{E:serie};
see~\cite{waldschmidt-adamczewski04}.
An example is the transcendence of
\begin{displaymath}
\sum_{n\ge0}g^{-2^n};
\end{displaymath}
see also \eqref{E:Kempner} above.

\begin{theorem}[Ridout, 1957]\label{T:Ridout}
For any $g\ge2$, any real algebraic number $\alpha$ and any $\epsilon>0$, the set of $p/q\in\bfQ$ with $q=g^k$ and $\vert\alpha-p/q\vert<q^{-1-\epsilon}$ is finite.
\end{theorem}

The theorems of Thue-Siegel-Roth and Ridout are very special cases of Schmidt's subspace
theorem~\cite{waldschmidt-schmidt80}, established in 1972, together with its $p$-adic extension by Schlickewei in 1976.
We state only a simplified version; see Bilu's Bourbaki lecture \cite{waldschmidt-bilu06} for references and further recent achievements based on this fundamental result.

For $\bfx=(x_0,\ldots,x_{m-1})\in\bfZ^m$, let $\vert\bfx\vert=\max\{\vert x_0\vert,\ldots,\vert x_{m-1}\vert\}$.

\begin{theorem}[Schmidt's subspace theorem]\label{T:SST}
Let $m\ge2$ be an integer, $S$ a finite set of places of $\bfQ$ containing the infinite place.
For each $v\in S$, let $L_{0,v},\ldots,L_{m-1,v}$ be $m$ independent linear forms in $m$ variables with algebraic coefficients in the completion of $\bfQ$ at~$v$.
Let $\epsilon>0$.
Then the set of $\bfx=(x_0,\ldots,x_{m-1})\in\bfZ^m$ for which
\begin{displaymath}
\prod_{v\in S}\vert L_{0,v}(\bfx)\ldots L_{m-1,v}(\bfx)\vert_v\le\vert\bfx\vert^{-\epsilon}
\end{displaymath}
is contained in the union of finitely many proper subspaces of~$\bfQ^m$.
\end{theorem}
    
Theorem~\ref{T:Roth} due to Thue-Siegel-Roth follows from Theorem~\ref{T:SST} if one takes
\begin{displaymath}
S=\{\infty\},
\hspace{.15in}
m=2,
\hspace{.15in}
L_0(x_0,x_1)=x_0,
\hspace{.15in}
L_1(x_0,x_1)=\alpha x_0-x_1.
\end{displaymath}
Theorem~\ref{T:Ridout} due to Ridout is also a consequence of Theorem~\ref{T:SST} if one takes
\begin{displaymath}
S=\{\infty\}\cup\{\ell:\ell\text{ prime and }\ell\mid g\}
\end{displaymath}
and
\begin{align}
L_{0,\infty}(x_0,x_1)=x_0,
\hspace{.15in}
&
L_{1,\infty}(x_0,x_1)=\alpha x_0-x_1,
\nonumber
\\
L_{0,\ell}(x_0,x_1)=x_0,
\hspace{.15in}
&
L_{1,\ell }(x_0,x_1)=x_1.
\nonumber
\end{align}
Since there is no loss of generality to assume that $g$ is squarefree, so that
\begin{displaymath}
g=\prod_{\ell\mid g}\ell,
\end{displaymath}
for $(x_0,x_1)=(q,p)$ with $q=g^k$, we have
\begin{align}
\vert L_{0,\infty}(x_0,x_1)\vert_{\infty}=q,
\hspace{.15in}
&
\vert L_{1,\infty}(x_0,x_1)\vert_{\infty}=\vert q\alpha-p\vert,
\nonumber
\\
\prod_{\ell\mid g}\vert L_{0,\ell}(x_0,x_1)\vert_\ell=q^{-1},
\hspace{.15in}
&
\prod_{\ell\mid g}\vert L_{1,\ell}(x_0,x_1)\vert_\ell=\vert p\vert_\ell\le1.
\nonumber
\end{align}

Further applications of the Subspace theorem to transcendence questions have been obtained by Corvaja and Zannier~\cite{waldschmidt-CZ02}.

%%%%%%%%%%%%%%%%%%%%
%
% SUBSECTION 4.2
%
%%%%%%%%%%%%%%%%%%%%

\subsection{Irrationality and Transcendence Measures}\label{SS:Irrationalitytranscendencemeasures}

The previous results can be made effective in order to reach irrationality measures or transcendence measures for automatic numbers.

In 2006, Adamczewski and Cassaigne \cite{waldschmidt-AC06} solved a conjecture of Shallit in 1999 by proving that the sequence of $g$-ary digits of a Liouville number cannot be generated by a finite automaton.
They obtained irrationality measures for automatic numbers.
Recall that the \emph{irrationality exponent} of an irrational real number $x$ is the least upper bound of the set of  numbers $\kappa$ for which the inequality
\begin{displaymath}
\left\vert x-\frac{p}{q}\right\vert<\frac{1}{q^\kappa}
\end{displaymath}
has infinitely many solutions $p/q$.

For instance, a Liouville number is a number whose irrationality exponent is infinite, while a real irrational algebraic number has irrationality exponent $2$ by Theorem~\ref{T:Roth}, as do almost all real numbers.

An explicit upper bound for the irrationality exponent for automatic irrational numbers is given by
\cite[Theorem~2.2]{waldschmidt-AC06}.
For the Prouhet-Thue-Morse-Mahler numbers for instance, the exponent of irrationality is at most~$5$.
However there is no uniform upper bound for such exponents, as pointed out to the author by Adamczewski.
The irrationality exponent for the automatic number associated with $\sigma(0)=0^n1$ and
$\sigma(1)=1^n0$ is at least~$n$.

Recently Adamczewski and Bugeaud \cite{waldschmidt-AB-a} have obtained transcendence measures for automatic numbers.
They show that automatic irrational numbers are either $S$- or $T$-numbers in Mahler's classification of transcendental numbers; see~\cite{waldschmidt-bugeaud04}.
This is a partial answer to a conjecture of Becker \cite{waldschmidt-AC06} which states that all automatic irrational numbers are $S$-numbers.

%%%%%%%%%%%%%%%%%%%%
%
% SECTION 5
%
%%%%%%%%%%%%%%%%%%%%

\section{Continued Fractions}\label{S:CF}

We discussed above some Diophantine problems which are related with the $g$-ary expansion of a real number.
Similar questions arise with the continued fraction expansion of real numbers.

%%%%%%%%%%%%%%%%%%%%
%
% SUBSECTION 5.1
%
%%%%%%%%%%%%%%%%%%%%

\subsection{Complexity of the Continued Fraction Expansion of an Algebraic Number}

In 1949, Khinchin asked the following question:
Are the partial quotients of the continued fraction expansion of a real non-quadratic irrational algebraic number bounded?
So far no example is known.
It is not yet ruled out that all these partial quotients are bounded for all such~$x$, or that they are unbounded for all such~$x$.
The common expectation seems to be that they are never bounded.
The situation in finite characteristic is quite different.
In 1976, Baum and Sweet \cite{waldschmidt-BS76} constructed a formal series which is cubic over the field $\bfF_2(X)$, the continued fractions of which have partial quotients of bounded degree.

We give here only a very short historical account.
The very first transcendence results for numbers given by their continued fraction expansions are due to  Liouville in 1844.
This topic was extensively developed by Maillet in 1906 and later by Perron in 1929.
In 1955, the Thue-Siegel-Roth theorem enabled Davenport and Roth to obtain deeper results.
Further investigations are due to Baker in 1962 and 1964.
The approximation results on real numbers by quadratic numbers due to Schmidt in 1967 are a main tool for the next steps by Davison in 1989, and by Queff\'elec \cite{waldschmidt-queffelec98} who established in 1998 the transcendence of the Prouhet-Thue-Morse continued fraction; see also
\cite[Section~13.7]{waldschmidt-AS03}.
There are further papers by Liardet and Stambul in 2000 and by Baxa in 2004.
In 2001, Allouche, Davison, Queff\'elec and Zamboni \cite{waldschmidt-ADQZ01} proved the transcendence of Sturmian continued fractions, namely continued fractions whose sequence of partial quotients is Sturmian.
The transcendence of the Rudin-Shapiro and of the Baum-Sweet  continued fractions was proved in 2005 by Adamczewski, Bugeaud and Davison.
In 2005, Adamczewski and Bugeaud \cite{waldschmidt-AB05a} showed that the continued fraction expansion of an algebraic number of degree at least three cannot be generated by a binary morphism.

%%%%%%%%%%%%%%%%%%%%
%
% SUBSECTION 5.2
%
%%%%%%%%%%%%%%%%%%%%

\subsection{The Fibonacci Continued Fraction}

The Fibonacci word discussed in Section~\ref{SS:FW} enabled Roy \cite{waldschmidt-roy03} to construct transcendental real numbers $\xi$ for which $\xi$ and $\xi^2$ are surprisingly well simultaneously and uniformly approximated by rational numbers.

Recall once more that $\Phi$ denotes the golden number, so that
\begin{displaymath}
\Phi^{-1}
=\Phi-1
=\frac{\sqrt{5}-1}{2}
=0.618\ldots.
\end{displaymath}

\begin{theorem}[Roy, 2003]\label{T:DRoy}
Let $A$ and $B$ be two distinct positive integers.
Let $\xi\in(0,1)$ be the real number whose continued fraction expansion is obtained from the Fibonacci word $w$ by replacing the letters $a$ and $b$ by $A$ and $B$ respectively to obtain
\begin{align}
&
[0;A,B,A,A,B,A,B,A,A,B,A,A,B,A,B,A,A,B,A,B,A,
\nonumber
\\
&
\hspace{.3in}
A,B,A,A,B,A,B,A,A,B,A,A,B\ldots].
\nonumber
\end{align}
Then there exists $c>0$ such that the inequalities
\begin{displaymath}
0<x_0\le X,
\hspace{.15in}
\vert x_0\xi-x_1\vert\le cX^{-\Phi^{-1}},
\hspace{.15in}
\vert x_0\xi^2-x_2\vert\le cX^{-\Phi^{-1}}
\end{displaymath}
have a solution in $\bfZ^3$ for all sufficiently large values of~$X$.
\end{theorem}

The arguments of Roy provide simplified proofs of Queff\'elec's results on the transcendence of the Thue-Morse and Fibonacci continued fractions.

Using Theorem~\ref{T:DRoy}, together with ideas of Davenport and Schmidt involving a transference theorem concerning Mahler's convex bodies, Roy also produced transcendental numbers which are surprisingly badly approximated by cubic algebraic integers.

Further results on Diophantine approximation of Sturmian continued fractions and on the simultaneous approximation of a number and its square have been obtained by Bugeaud and Laurent, Fischler, Roy, and more recently by Adamczewski and Bugeaud~\cite{waldschmidt-AB-a}.

For further references on Diophantine approximation, we refer the reader to Bugeaud's recent
book~\cite{waldschmidt-bugeaud04}.

%%%%%%%%%%%%%%%%%%%%
%
% SECTION 6
%
%%%%%%%%%%%%%%%%%%%%

\section{Open problems}\label{S:OP}

Many problems related to the present topic are open.
To trace their original sources would require more thorough bibliographical investigation.
While the origins of these questions are most often much older, we give references to recent papers where they are quoted.

Among the open questions raised in \cite{waldschmidt-AB05a} is the following:
\begin{itemize}
\item Does there exist an algebraic number of degree at least $3$ whose continued fraction expansion is generated by a morphism?
\end{itemize}
In the same vein the next question is open:
\begin{itemize}
\item Do there exist an integer $g\ge3$ and an algebraic number of degree at least $3$ whose   expansion in base $g$ is generated by a morphism?
\end{itemize}

Two open questions, for which positive answers are expected, are proposed
in~\cite{waldschmidt-AR-a}:
\begin{itemize}
\item Is it true that the binary expansion of every algebraic number contains arbitrarily large squares?
\item Is it true that the binary expansion of every algebraic number contains arbitrarily large palindromes?
\end{itemize}
Recall that a square is nothing other than a $2$-power, namely a pattern $XX$ where $X$ is a word, while a palindrome is a word $W=w_1w_2\ldots w_r$ which is invariant under reversal: if $\overline{W}$ denotes the word $w_r\ldots w_2w_1$, then $W=\overline{W}$.
Hence a palindrome is either of the form $X\overline{X}$ or of the form  $Xx\overline{X}$ where $X$ is a word and $x$ a letter.

Other problems are suggested by Rivoal~\cite{waldschmidt-rivoal-a}:
\begin{itemize}
\item Let $g\ge2$ be an integer.
Give an explicit example of a real number $x>0$ which is simply normal in base $g$ and for which $1/x$ is not simply normal in base~$g$.
\item Let $g\ge2$ be an integer.
Give an explicit example of a real number $x>0$ which is normal in base $g$ and for which $1/x$ is not normal in base~$g$.
\item Give an explicit example of a real number $x>0$ which is normal and for which $1/x$ is not normal.
\end{itemize}

\begin{remark}
In~\cite{waldschmidt-LST96}, there is a construction of an automatic number, the inverse of which is not automatic.
This answers by anticipation \cite[Section~13.9, Problem~2]{waldschmidt-AS03}.
\end{remark}

From the open problems in \cite[Section~13.9]{waldschmidt-AS03}, we select the following two:
\begin{itemize}
\item Show that the number $\log 2$ is not $2$-automatic; see \eqref{E:log2}.
\item Show that the number $\pi$ is not $2$-automatic; see \eqref{E:pi}.
\end{itemize}

We conclude this paper with one last open problem, attributed to Mahler; see for
instance~\cite{waldschmidt-AB06}:
\begin{itemize}
\item Let $(e_n )_{n\ge1}$ be an infinite sequence over $\{0,1\}$ that is not ultimately periodic.
Is it true that at least one of the two numbers
\begin{displaymath}
\sum_{n\ge1}e_n2^{-n}
\hspace{.15in}
\text{and}
\hspace{.15in}
\sum_{n\ge1}e_n3^{-n}
\end{displaymath}
is transcendental?
\end{itemize}
From Conjecture \ref{Conj:1} with $g=3$ and $a=2$, it follows that the second number should be always  transcendental.

\begin{acknowledgments}
The author is grateful to Boris Adamczewski, Jean-Paul Allouche, Yann Bugeaud, St\'ephane Fischler, Damien Roy, Tanguy Rivoal and Paul Voutier for their useful comments on a preliminary version of this paper.
This topic was the subject of a lecture by the author at a conference in honour of the 65th birthday of 
Professor Iekata Shiokawa, held at Keio University, Yokohama, Japan, on 7--10 March 2006.
The author is also grateful to Takao Komatsu for his invitation to participate at this conference. 
A tremendous job was done by William Chen who edited and improved substantially the submitted paper before publication in this volume. 
\end{acknowledgments}

\end{document}